\documentstyle[epsf,12pt]{article}
\newcommand{\beq}{\begin{equation}}
\newcommand{\eeq}{\end{equation}}
\title{On the Accurate Finite Element Solution of a Class of Fourth
Order Eigenvalue Problems}
\author{B.M.\ Brown\thanks{Department of Computer Science, Cardiff
University of Wales, Cardiff CF2 3XF, Wales, UK.}
\hspace*{1em}
E.B.\ Davies\thanks{Department of Mathematics, King's College London,
Strand, London WC2R 2LS, UK.}
\hspace*{1em}
P.K.\ Jimack\thanks{School of Computer Studies,
University of Leeds, Leeds LS2 9JT, UK.}
\hspace*{1em}
M.D.\ Mihajlovi\'c\footnotemark[1]}
\date{\vspace{-11pt}}
\begin{document}
\maketitle

\begin{abstract}
This paper is concerned with the accurate numerical approximation of the
spectral properties of the biharmonic operator on various domains in two
dimensions. A number of analytic results concerning the eigenfunctions of
this operator are summarized and their implications for numerical
approximation are discussed. In particular, the asymptotic behaviour of the
first eigenfunction is studied since it is known that this has an unbounded
number of oscillations when approaching certain types of corner on domain
boundaries. Recent computational results of Bj\o rstad and Tj\o stheim
\cite{Bjo}, using a highly accurate
spectral Legendre-Galerkin method, have demonstrated that a number of these
sign changes may be accurately computed on a square domain {\em provided}
sufficient care is taken with the numerical method. We demonstrate that
similar accuracy is also achieved using an unstructured finite element
solver which may be applied to problems on domains with arbitrary
geometries. A number of results obtained from this mixed finite element
approach are then presented for a variety of domains. These include
a family of circular sector regions, for which the oscillatory behaviour
is studied as a function of the internal angle, and another family
of (symmetric and non-convex) domains, for which the parity of the least 
eigenfunction
is investigated. The paper not only verifies existing asymptotic theory,
but also allows us to make a new conjecture concerning  the eigenfunctions
of the biharmonic operator.
\end{abstract}

\section{Introduction}
In recent years there has been growing interest  in the
 spectral theory of higher order
elliptic operators. However, in contrast to the  well-understood
theory for second order operators, the higher order theory
can be quite different and is certainly less well developed.
The recent survey article \cite{Dav} contains an  account of many
of  the key results. The purpose of this paper is to explore, using
reliable and accurate numerical techniques, some of the properties of the
eigenfunctions of the biharmonic operator
on domains in $\Re^2$. In particular we investigate
the existence of so-called {\em nodal lines}
in the neighbourhood of certain corners of the domain and also how the
 parity of these eigenfunctions may change with domain geometry.
As will become clear from the quantitative description  of these problems below
this is a demanding computational task since very high numerical accuracy
is required in order to resolve the phenomena under investigation.

The biharmonic eigenvalue problem typically comes in two different forms:
the clamped plate eigenproblem
\begin{equation}
\Delta^2 u = \lambda u \;, \label{eq:1.1}
\end{equation}
and the buckling plate eigenproblem
\begin{equation}
\Delta^2 u = \lambda \Delta u \;, \label{eq:1.2}
\end{equation}
each on a domain $\Omega \subset \Re^2$ subject to appropriate conditions on
the boundary $\partial \Omega$. Throughout this paper we shall
consider only the zero Dirichlet conditions
\begin{equation}
u=\frac{\partial u}{\partial n}=0 \;\;\;\;\forall x \in \partial \Omega.
\end{equation}
\par
In the remainder of this introductory section we describe details of
 the specific problems that are to be considered.
Section \ref{s:2}  then contains  an outline and justification of the numerical 
techniques,
based upon the use of conforming $C^0$ finite elements, that are used to undertake
these investigations. This is followed in Sections \ref{s:3} and \ref{s:4} 
respectively
by our findings concerning oscillations of the eigenfunctions in the neighbourhood
of  corners of the domain and the dependence of the parity of the eigenfunctions 
on
the domain geometry.
The paper concludes with a brief discussion of our results along with a 
consideration
of the issues associated with the validation of  numerical simulations such as 
these,
including mesh convergence and the use of rigorous enclosure~techniques.

\subsection{Oscillatory properties of eigenfunctions} \label{s:1.1}
For each of the two problems (\ref{eq:1.1}) and (\ref{eq:1.2})
it is shown in \cite{B+R},\cite{Coff} that in the neighbourhood of a corner
of $\partial \Omega$ with sufficiently small internal angle $\theta$
any eigenfunction changes sign an
infinite number of times. Moreover the ratio of the distance from the corner, 
along
its bisector, of consecutive zeros of the eigenfunction tends to a limit which
is dependent upon $\theta$. It transpires that this ratio is an increasing 
function of $\theta$ and is such that the oscillations cease to be present
when the internal angle exceeds some critical value $\theta_c$.

These results are  special cases of more general theorems concerning higher
order elliptic operators in greater than one dimension \cite{KKM}. Their
derivation is based
upon an asymptotic expansion of the eigenfunction, $u$, centred at, and
in a neighbourhood of, the corner. The leading term of such an expansion
as one approaches the corner along the bisector of the angle (which we
always do)
is of the form $cr^p$, where $r$ is the distance from the corner,
and $p$ is a solution of the transcendental equation
\begin{equation}
p+1+\frac{\sin{((p+1)\theta )}}{\sin \theta }=0. \label{eq:1.4}
\end{equation}
Assuming that the coefficient $c$ is non-zero the relevant
solution is the one with smallest real part. From this equation,
which has only real solutions for
$\theta\!>\!\theta_c\!=\!0.8128 \pi\!=\!146^\circ 30^\prime$, it may be 
concluded that when the internal angle is greater than $\theta_c$ the 
eigenfunction has no sign changes near the corner. Conversely as 
$\theta\!\rightarrow\!0$
both the real and imaginary parts of the exponent of the
leading term in the asymptotic expansion grow unboundedly, implying
that the eigenfunction oscillates with high frequency but is damped out
quickly for small $\theta$.
Moreover, when we let $s_n$ be the distance along the bisector
from the corner to the $n^{th}$ zero of $u$, where $n$ increases with decreasing 
distance
from the corner, then it follows immediately from the asymptotic
expression that
\begin{equation}
\frac{s_n}{s_{n+1}} \sim e^{\pi/\beta} \;\;\;\; \mbox{as $n \rightarrow \infty$},
\label{eq:1.5}
\end{equation}
where $p=\alpha+i \beta$ is the solution of (\ref{eq:1.4}) with smallest real 
part.
Furthermore, when $r_n$ is the distance along the bisector to the $n^{th}$ 
extremum  of $u$
and $t_n$ is the magnitude of this extremum, it may also be shown that
\begin{equation}
\frac{r_n}{r_{n+1}} \sim \frac{s_n}{s_{n+1}} \;\;\;\; \mbox{ and }
\;\;\;\;\frac{t_n}{t_{n+1}} \sim \left( \frac{s_n}{s_{n+1}} \right)^\alpha
\label{eq:1.6}
\end{equation}
as $n \rightarrow \infty$.
The first goal of this research is to verify this behaviour
for a sequence of domains which take the form of sectors of a circle
with increasing arc length. In view of the above discussion and the data
in Table~\ref{t0} below it is clear that such a verification will require
very high accuracy from the numerical methods used. This table presents
the solutions $p=\alpha+ i\beta$ of (\ref{eq:1.4}) as a function of
$\theta$, as well as the ratios $e^{\pi/\beta}$ and $e^{\alpha\pi/\beta}$
which appear in (\ref{eq:1.5}) and (\ref{eq:1.6}). These solutions of
(\ref{eq:1.4}) are obtained by   Newton's method with double precision
arithmetic.

\begin{table}[htb]
\begin{center}
{\small\begin{tabular}{|c|c|c|c|c|c|c|c|}\hline
$\theta$ & $10^\circ$ & $20^\circ$ & $30^\circ$ & $40^\circ$ & $50^\circ$ &
$60^\circ$ & $70^\circ$ \\ \hline
$\alpha=\Re{\rm e}(p)$ & 25.141144 & 13.079480 & 9.062965 & 7.057831 &
5.857356 & 5.059329 & 4.491404 \\ \hline
$\beta=\Im{\rm m}(p)$ & 12.864086 & 6.384388 & 4.202867 & 3.095366 &
2.416840 & 1.952050 & 1.608491 \\ \hline
$e^{\pi/\beta}$ & 1.27662 & 1.63571 & 2.11169 & 2.75918 & 3.66884 & 4.99972 &
7.05073 \\ \hline
$e^{\alpha\pi/\beta}$ & 463.97239 & 623.95258 & 875.20459 & 1291.07923 &
2026.03589 & 3437.12115 & 6452.99420 \\ \hline
\end{tabular}}

\medskip

{\small\begin{tabular}{|c|c|c|c|c|c|c|c|}\hline
$\theta$ & $80^\circ$ & $90^\circ$ & $100^\circ$ & $110^\circ$ & $120^\circ$ &
$130^\circ$ & $140^\circ$ \\ \hline
$\alpha=\Re{\rm e}(p)$ & 4.067435 & 3.739593 & 3.479215 & 3.268096 & 3.094139 &
2.949023 & 2.826869 \\ \hline
$\beta=\Im{\rm m}(p)$ & 1.339586 & 1.119024 & 0.930373 & 0.762118 & 0.604585 &
0.446356 & 0.261695 \\ \hline
$e^{\pi/\beta}$ & 10.43532 & 16.56743 & 29.27404 & 61.69387 & 180.5992 & 1139.464 &
163533.23 \\ \hline
$e^{\alpha\pi/\beta}$ & 13890.112 & 36267.559 & 126534.61 & 709058.99 &
9607097.6 & 1033431725 & $0.5472\cdot10^{15}$\\ \hline
\end{tabular}}
\caption{The solutions of (\ref{eq:1.4}) as a function of
the angle $\theta$. The ratios $e^{\pi/\beta}$ and $e^{\alpha\pi/\beta}$
are also given. \label{t0}}
\end{center}
\end{table}

\subsection{Parity of eigenfunctions for non-convex domains} \label{s:1.2}
Our other main goal is to investigate the dependence of the parity of the
eigenfunctions  of problems (\ref{eq:1.1}) and (\ref{eq:1.2}) on the geometry
of $\Omega$.  Given a domain $\Omega \subset \Re^2$ which is symmetric about
some line (which, without loss of generality, we will choose to be the
$y$-axis), it is known that the least eigenfunction of the {\em Laplace}
operator, subject to zero Dirichlet boundary conditions, is always positive and
of multiplicity one \cite{dav2}. Symmetry arguments now force it to be an even
function of $x$.  No equivalent result holds for the biharmonic operator.
Indeed, it is the purpose of this work to demonstrate numerically that
no such result holds for problems of the form (\ref{eq:1.1}) and (\ref{eq:1.2})
in arbitrary symmetric domains $\Omega$. This is achieved in  Section \ref{s:4}
by considering a family of non-convex domains whose boundary is defined by
the closed curves
\begin{equation}
 y=\pm ( c + x^2-x^4) \label{eq:1.7}
\end{equation}
for $ c  >0$ (see figure \ref{fig:1} below).
\begin{figure}[htbp]
\centerline{
\epsfxsize=11cm
\epsffile{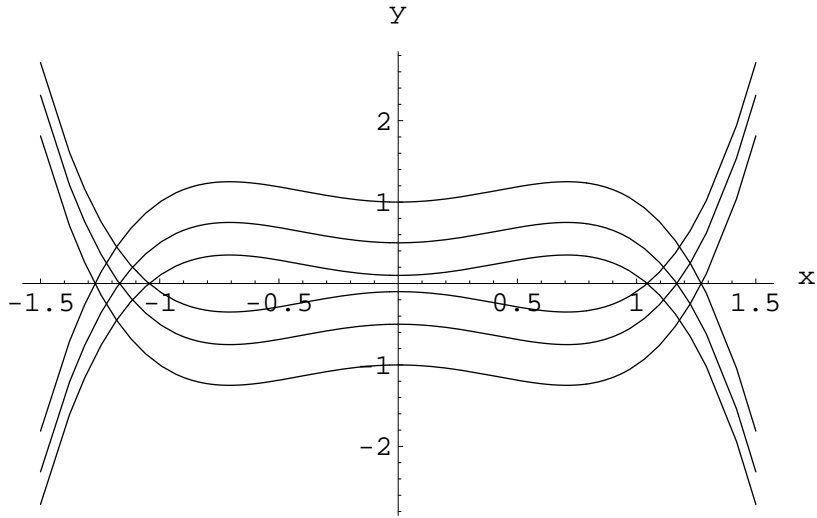}
}
\caption{The dumb-bell shape domain for $c=1,\; c=0.5,\; c=0.1$}
\label{fig:1}
\end{figure}
 Domain monotonicity arguments \cite{dav2} imply that
the eigenvalues of the biharmonic operator increase as $c\to 0$ and converge to 
the
eigenvalues of the same operator on the region associated with $c=0$.
This region is disconnected and its
eigenvalues have multiplicity two. Variational
considerations now show that for small enough $c>0$ the low-lying
eigenvalues occur in pairs, one associated with an even
eigenfunction and the other associated with an odd eigenfunction.
However in the biharmonic case one cannot say which of these is
the smaller on purely variational grounds.
Indeed by monitoring the smallest two eigenvalues of (\ref{eq:1.1}) or 
(\ref{eq:1.2})
as functions of $c$, it may be observed that their paths cross on more
than one occasion as $c \rightarrow 0$. Furthermore it is apparent
that one of these eigenvalues corresponds to an even eigenfunction
while the other corresponds to an odd one, and consequently we observe
that the parity of the least eigenfunction is not fixed as a function of $c$.
We again remark on the need for very high accuracy in our computational
algorithms; this time in order to identify values of $c$ for which the
eigenvalue paths cross.
\section{Numerical method} \label{s:2}
Given the above discussion, it is clear that our discretization scheme
must be able to accommodate quite general domain geometries. For
this reason we choose a numerical scheme which is based upon the use of
conforming $C^0$ Lagrange finite elements (see \cite{Str},\cite{Zien}
for example) on some triangulation of the domain $\Omega$. The particular
method that we use is based on that of Ciarlet and Raviart in~\cite{C-R}.

Let $H^1(\Omega)$ represent the usual Sobolev space of $L^2(\Omega)$
functions whose first partial derivatives are in $L^2(\Omega)$, and
\[
H^1_0(\Omega)=\left\{\psi\in H^1(\Omega):\psi|_{\partial\Omega}=
\left.\frac{\partial \psi}{\partial n}\right|_{\partial\Omega}=0\right\}.
\]
Consider a family of regular and quasi-uniform triangulations,
${\cal T}^h$ (where $h$ is the diameter of the largest triangle),
of the domain $\Omega$ (see, for example, Ciarlet \cite{Cia}), then we may
define the space
\[
{\cal S}^h_m=\{p\in C^0(\bar\Omega)\,:\, p|_T\in P_m(T),\
\forall T\in{\cal T}^h\},
\]
 where $P_m(T)$ is the space of polynomials of degree at most $m$ over
triangle $T$. From this we may define the following finite-dimensional
subspaces of $H^1(\Omega)$: $V^h={\cal S}^h_m$ and
$W^h={\cal S}^h_m\cap H^1_0(\Omega)$. Then the Ciarlet--Raviart mixed
finite element method approximates the solution $u$ of (\ref{eq:1.1})
by the component $u^h$ of  the solution $\{v^h,u^h\} \in V^h \times W^h$
to the following problem.
\newtheorem{t2.1}{Problem}[section]
\begin{t2.1}
Find $\{v^h,u^h\} \in V^h \times W^h$ such that
\begin{eqnarray*}
\int_\Omega v^h w_1 \;d \Omega + \int_\Omega \nabla u^h \cdot \nabla w_1 \; 
d\Omega
  & = & 0 \;\;\;\;\qquad\ \forall w_1\in V^h,\\
\int_\Omega \nabla v^h \cdot \nabla w_2 \; d\Omega
  & = & -\lambda\int_\Omega u^h w_2 \; d \Omega \;\;\;\; \forall w_2\in W^h.
\end{eqnarray*} \label{p:2.1}
\end{t2.1}

In this work we choose $m=3$, i.e.\ piecewise cubic approximation, since such a
discretization is known to  satisfy the Brezzi  stability condition, \cite{Bre}, 
\cite{Shaid} and has
the following asymptotic error property (\cite{C-R}):
\begin{equation}
\| u-u^h \|_{H^1(\Omega)} + \| \Delta u - v^h \|_{L_2(\Omega)}
\leq c \| u \|_{H^{5}(\Omega)} h^{2}\;,
\label{err-est}
\end{equation}
provided the solution $u$ is sufficiently smooth. Hence the trial
functions $v^h$ and $u^h$ in Problem \ref{p:2.1} may be expressed as
\begin{equation}
v^h=\sum^{n_i+n_d}_{j=1}v_j\phi_j \qquad\mbox{and}\qquad
u^h=\sum^{n_i}_{j=1}u_j\phi_j, \label{e:trial}
\end{equation}
where the $\phi_j$'s are the usual piecewise cubic Lagrange basis functions
on ${\cal{T}}^h$ (see \cite{Str},\cite{Zien} for example), and $n_i$ and
$n_d$ are the number of nodes in the interior of $\Omega$
and on its Dirichlet boundary respectively. (We also use the convention that
the interior nodes are numbered prior to those on the Dirichlet boundary.)
This discretization leads to the following matrix problem
\begin{equation}
\left[\begin{array}{cc} M&K^t\\ K&O\end{array}\right]\left[\begin{array}{c}
\underline v\\ \underline u\end{array}\right]=- \lambda \left[\begin{array}{cc}
\underline 0\\ \hat{M} \underline u \end{array}\right],
\label{e2.7}
\end{equation}
where $M \in \Re^{(n_i + n_d) \times (n_i + n_d)}$,
$\hat{M} \in \Re^{ n_i \times n_i}$ is a principal submatrix of $M$, $K \in 
\Re^{n_i \times ( n_i + n_d) }$ and $K^t$ is the transpose of $K$. The specific
entries of these matrices are given by
\begin{eqnarray*}
 K_{ij} &=&\int_\Omega \nabla \phi_i \cdot \nabla \phi_j \; d \Omega, \\
 M_{ij} &=&\int_\Omega \phi_i \phi_j \; d\Omega, 
\end{eqnarray*}
and the vectors
$\underline v \in \Re^{n_i+n_d}$ and $\underline u \in \Re^{n_i}$ denote
the coefficients of $v^h$ and $u^h$ respectively in (\ref{e:trial}).

It may be observed that $\underline{u}$ and $\lambda$ in (\ref{e2.7})
are generalized eigenvectors and eigenvalues respectively for the corresponding
Schur complement problem:
\begin{equation}
KM^{-1}K^t \underline{u}= \lambda \hat{M}\underline{u}. \label{e:SCM}
\end{equation}
In this paper we are typically interested in approximating the smallest
eigenvalues of this problem and so make use of straightforward inverse
iteration (as described in \cite{Golub} for example).
Since the Schur complement matrix $KM^{-1}K^t$ is dense, whereas the matrix
on the left-hand side of (\ref{e2.7}) is sparse, it is convenient to apply
the inverse iteration to the latter system. This yields the following
matrix equation at the $(n+1)^{th}$  inverse iteration:
\begin{equation}
\left[\begin{array}{cc} M&K^t\\ K&O\end{array}\right]\left[\begin{array}{c}
\underline v^{n+1}\\ \underline u^{n+1}\end{array}\right]=  
\left[\begin{array}{cc}
\underline 0\\ -\hat{M} \underline u^n \end{array}\right].
\label{e2.8}
\end{equation}
When approximating the second smallest eigenpair, as in Section \ref{s:4} for 
example,
this iteration is slightly modified by subtracting out the component
of the smallest eigenvector from $\underline{u}^{n+1}$ at each
iteration.

We solve  the systems (\ref{e2.8}) using a direct sparse method
based upon an initial block reduction followed by a sparse Cholesky
factorization of symmetric positive definite sub-blocks as described
in \cite{BJM99}. This permits
re-use of the factorizations for the second and subsequent solves,
which ensures that the cost of these is significantly less than that
of the initial solution of (\ref{e2.8}). An additional advantage of a
direct, rather than an iterative, method such as this is that it is
possible to obtain highly accurate solutions: something that is required
in both Sections \ref{s:3} and \ref{s:4} below.  To demonstrate the
accuracy of this approach (both the discretization and the
direct solution), we now present some typical numerical computations
for a biharmonic eigenvalue problem which has been considered
in some detail elsewhere, \cite{Bjo},\cite{BJM99}.

Recall from Subsection \ref{s:1.1} that in the neighbourhood of a corner with
a sufficiently small internal angle, the eigenfunctions for both problems 
(\ref{eq:1.1}) and (\ref{eq:1.2}) are known to change sign at an infinite
number of points.  In particular, when the internal angle $\theta=\pi/2$
the solution of (\ref{eq:1.4}) with smallest real part is
$3.739593284 + 1.1190248i$ (see Table \ref{t0}). It then follows
from (\ref{eq:1.5}) and (\ref{eq:1.6}) that the asymptotic ratio
of distances along the bisector of consecutive zeros and extremal points of
the eigenfunctions should both tend to $16.567428$ as $n \rightarrow \infty$.
Also from (\ref{eq:1.6}) we see that the ratios between the magnitudes of
consecutive extrema along the bisector tends to $36267.550$ as $n \rightarrow 
\infty$.
In recent years many authors have tried to verify numerically these
theoretical results using a variety of computational techniques (see, for
example, \cite{Beh},\cite{Chen},\cite{Hack},\cite{Wien}). Since the oscillatory
features occur very close to
the corners and are damped out very quickly, most of these attempts
have, due to discretization errors and numerical inaccuracy, failed to find
more than one sign change. A notable exception to this is the recent paper
of Bj\o rstad
and Tj\o stheim (\cite{Bjo}) in which the authors report five correct sign
changes for the principal eigenfunction of problem (\ref{eq:1.1}) on the
 domain $(0,1)\times (0,1)$. This is achieved using a spectral
Legendre--Galerkin method  and by \mbox{performing computations with
quadruple precision arithmetic.}

In Tables \ref{t:1} to \ref{t:4} we present some comparative
numerical results of our own, obtained using the mixed finite element
method described above on three different meshes. In addition we contrast in
the last column of the tables \ref{t:1}, \ref{t:3} and \ref{t:4} the best
results of Bj\o rstad and Tj\o stheim reported in \cite{Bjo}. (This
latter paper does not present results for the positions, $s_n$, of the
zeros and so no comparisons are possible for these figures). It may clearly 
be seen that, for each of these meshes, we are also able to identify at least
five sign changes for the principal eigenfunction of (\ref{eq:1.1})
on the unit square. These calculations have been undertaken using double
precision arithmetic on one eighth of the domain, using symmetries at
$x=1/2$, $y=1/2$ and $y=x$, and with heavy local refinement of the meshes
near the corner ($h\approx 10^{-7}$ in the vicinity of the corner
with a gradual transition to $h\approx 10^{-2}$ at the centre of the domain).
The results   are  entirely consistent with the quadruple precision
results presented in \cite{Bjo} and with the predicted asymptotic ratios
(\ref{eq:1.5}) and (\ref{eq:1.6}). This leads us to have a reasonable degree of
confidence in the underlying numerical procedure that we use for the 
investigations undertaken in the following two sections.

\begin{table}[htb]
\begin{center}
{\footnotesize \begin{tabular}{|c|c|c|c|c|} \hline
& \multicolumn{3}{c|}{FE method} & L--G method \cite{Bjo} \\ \hline
$N$ & 5575 & 22351 & 90028 & 5000\\ \hline
$\lambda_1$ & 1294.934809196 & 1294.934021432 & 1294.933981245 & 1294.933979592
\\ \hline \hline
$s_1$ & 0.042311747881 & 0.042310932855 & 0.042310963855 & --- \\ \hline
$r_1$ & 0.032634299694 & 0.032630435492 & 0.032629530244 & 0.032629472978\\ \hline
$t_1$ & $-0.169815505603\cdot10^{-4\phantom{0}}$ & $-0.169795463221\cdot10^{-4
\phantom{0}}$ & $-0.169791420686\cdot10^{-4\phantom{0}}$ & 
$-0.169791287981\cdot10^{-4\phantom{0}}$ \\ \hline \hline
$s_2$ & 0.002553961304 & 0.002553862560 & 0.002553860600 & --- \\ \hline
$r_2$ & 0.001969665622 & 0.001969564641 & 0.001969500077 & 0.001969491936 \\ \hline
$t_2$ & $\ 0.468412256179\cdot10^{-9\phantom{0}}$ & $\ 0.468173039472\cdot10^{-9
\phantom{0}}$ & $\ 0.468161662361\cdot10^{-9\phantom{0}}$ &
$\ 0.468161006275\cdot10^{-9\phantom{0}}$ \\ \hline \hline
$s_3$ & 0.000154151151 & 0.000154149893 & 0.000154149497 & --- \\ \hline
$r_3$ & 0.000118849965 & 0.000118881375 & 0.000118877347 & 0.000118877352\\ \hline
$t_3$ & $-0.129136601165\cdot10^{-13}$ & $-0.129088794006\cdot10^{-13}$ &
$-0.129085648295\cdot10^{-13}$ & $-0.129085369432\cdot10^{-13}$\\ \hline \hline
$s_4$ & 0.000009304900 & 0.000009304402 & 0.000009304373 & ---\\ \hline
$r_4$ & 0.000007177212 & 0.000007175502 & 0.000007175314 & 0.000007175365 \\ \hline
$t_4$ & $\ 0.356091995705\cdot10^{-18}$ & $\ 0.355934162890\cdot 10^{-18}$ &
$\ 0.355925959506\cdot 10^{-18}$ & $\ 0.355925258182\cdot10^{-18}$\\ \hline \hline
$s_5$ & 0.000000561678 & 0.000000561628 & 0.000000561569 & --- \\ \hline
$r_5$ & 0.000000433175 & 0.000000432837 & 0.000000432811 & 0.000000432767\\ \hline
$t_5$ & $-0.982438718508\cdot 10^{-23}$ & $-0.981777297429\cdot10^{-23}$ &
$-0.981303186480\cdot10^{-23}$ & $-0.974860961611\cdot10^{-23}$\\ \hline
\end{tabular}}
\caption{Positions of the local zeros ($s_n$) and local extrema
($r_n$), and the values of the local extrema ($t_n$), for the first
eigenfunction of the clamped plate problem (\ref{eq:1.1}) for $n=1$ to $5$,
as calculated on three different meshes (where $N$ is the number of
unknowns in the system (\ref{e2.7}) for each mesh). The last column
contains the corresponding results from [4] where available.
\label{t:1}}
\end{center}
\end{table}
\begin{table}[htb]
\begin{center}
{\small\begin{tabular}{|c|c|c|c|} \hline
$N$ & 5575 & 22351 & 90028 \\ \hline
$s_1/s_2$ & 16.5671 & 16.5674  & 16.5675 \\ \hline
$s_2/s_3$ & 16.5679 & 16.5674  & 16.5674 \\ \hline
$s_3/s_4$ & 16.5667 & 16.5674  & 16.5674 \\ \hline
$s_4/s_5$ & 16.5663 & 16.5668  & 16.5685 \\ \hline \hline
limit & 16.5674  & 16.5674  & 16.5674  \\ \hline
\end{tabular}}
\caption{Ratios between consecutive local zeros as computed on three
different meshes (where $N$ is the number of unknowns in the system
(\ref{e2.7}) for each mesh). \label{t:2}}
\end{center}
\end{table}\nopagebreak
\begin{table}[htb]
\begin{center}
{\small\begin{tabular}{|c|c|c|c|c|} \hline
& \multicolumn{3}{c|}{FE method} & L--G method \cite{Bjo} \\ \hline
$N$ & 5575 & 22351 & 90028 & 5000 \\ \hline
$r_1/r_2$ & 16.5684 & 16.5673 & 16.5674 & 16.5675 \\ \hline
$r_2/r_3$ & 16.5727 & 16.5675 & 16.5675 & 16.5674 \\ \hline
$r_3/r_4$ & 16.5593 & 16.5677 & 16.5675 & 16.5674\\ \hline
$r_4/r_5$ & 16.5689 & 16.5778 & 16.5784 & 16.5801 \\ \hline \hline
limit & 16.5674 & 16.5674 & 16.5674 & 16.5674\\ \hline
\end{tabular}}
\caption{Ratios between consecutive local extremal positions as computed
on three different meshes (where $N$ is the number of unknowns in the system
(\ref{e2.7}) and the system from [4] respectively).
\label{t:3}}
\end{center}
\end{table}
\begin{table}[htb]
\begin{center}
{\small\begin{tabular}{|c|c|c|c|c|} \hline
& \multicolumn{3}{c|}{FE method} & L--G method \cite{Bjo} \\ \hline
$N$ & 5575 & 22351 & 90028 & 5000 \\ \hline
$|t_1/t_2|$ & 36253.4292 & 36267.6722 & 36267.6901 & 36267.7125 \\ \hline
$|t_2/t_3|$ & 36272.6177 & 36267.5198 & 36267.5223 & 36267.5498\\ \hline
$|t_3/t_4|$ & 36264.9548 & 36267.6044 & 36267.5564 & 36267.5496\\ \hline
$|t_4/t_5|$ & 36245.7209 & 36254.0633 & 36270.7433 & 36510.3612\\ \hline \hline
limit & 36267.5596  & 36267.5596 & 36267.5596 & 36267.5596 \\ \hline
\end{tabular}}
\caption{Ratios between consecutive local extremal values as computed on three 
different meshes (where $N$ is the number of unknowns in the system
(\ref{e2.7}) and the system from [4] respectively).
\label{t:4}}
\end{center}
\end{table}

\section{Eigenfunctions on a circular sector} \label{s:3}

The purpose of this section is to extend the numerical studies in
\cite{Bjo},\cite{BJM99},
which concern the behaviour of the principle eigenfunction
of the biharmonic operator near
the corners of a square domain, to the case of
a general internal angle $\theta$. This is achieved  by considering a one 
parameter
family of unit radius circular sector domains, $\Omega_\theta$, and applying
to (\ref{eq:1.1}) the finite
element discretization defined in Problem \ref{p:2.1}.
Recall from Subsection \ref{s:1.1} that the oscillatory behaviour  of
the eigenfunction, $u$, near to a corner depends upon the internal angle $\theta$
through the imaginary part of a solution of equation (\ref{eq:1.4}).
In particular, there is
a critical angle, $\theta_c$, above which no oscillations are
present. As usual
all of our calculations are restricted to the behaviour of the
eigenfunction on the bisector of the angle at the corner.
In this study we examine in detail the behaviour of the oscillations in the 
principal
eigenfunction as $\theta$ approaches $\theta_c=0.8128 \pi=146^\circ 30^\prime$,
but we also present approximation of the eigenfunction for other values of
the angle $\theta$.

It has already been demonstrated in Section \ref{s:2} that our finite
element discretization scheme is capable of producing the accuracy
that is required to resolve up to five sign changes
when $\theta\!=\!\pi/2$, on the unit square, which is consistent with the best
numerical results of which we are aware \cite{Bjo}.  Clearly a key to obtaining 
results
with this level of resolution near to the corner is the  use
of an appropriate triangulation ${\cal {T}}^h$.
For the circular sector regions that we consider here we use a similar
mesh generation strategy to that adopted for the one-eighth of the unit square
in the previous section. This involves the creation of an unstructured
mesh based upon polar coordinates centred at the vertex. When $r$ is greater
than some (very small) chosen value, $\rho_1$ say, the mesh is such that
$h \approx r h_1$, and when $h < \rho_1$ the mesh is approximately uniform. 
In this approach it is necessary to represent the circular arc by a
piecewise affine approximation of side length $\!\approx\!h_1$. In order
to obtain an idea of the magnitude of the errors resulting from such an
approximation we choose to perform all calculations on two different
polygonal domains: one containing the sector and the other contained by it.
It follows from the min-max principle that the $m^{th}$ eigenvalue of the
differential operator  on the sector is contained between the corresponding
eigenvalue on each of these computational domains.

In addition to considering the errors that result from taking an approximation
to the domain $\Omega_\theta$, it is also desirable to have an indication of the
errors associated with approximating the continuous operator by the finite
element discretization.
In the numerical results that follow the calculations for each different value
of $\theta$ are presented for two different pairs of meshes of the type
described above. Typically, the first pair of such meshes has only about $30\%$
of the number of degrees of freedom as the second. Hence, when intervals for
the local zeros/extrema positions 
obtained on the finer pair of meshes are contained within those obtained
on the coarser pair, we can have a reasonable expectation of
the reliability of these calculations.

Tables \ref{t:3.1} and \ref{t:3.2} show the computed values of
$s_n$ and $r_n$ ($n=1$ to $4$), the positions of the zeros and extrema
respectively of the first eigenfunction of (\ref{eq:1.1}), for domains
$\Omega_\theta$ with $\theta$ between $30^\circ$ and $140^\circ$ inclusive. 
 
\begin{table}[htb]
\begin{center}
{\small\begin{tabular}{|c|c|c|c|c|c|c|} \hline
$\alpha$ & $30^\circ$ & $40^\circ$ & $50^\circ$ & $60^\circ$ & $70^\circ$ &
$80^\circ$ \\ \hline
$h_1$ & $0.065$ & $0.087$ & $0.110$ & $0.132$ & $0.154$ & $0.176$ \\ \hline
$\rho_1$ & $5\cdot10^{-4}$ & $5\cdot10^{-5}$ & $1\cdot10^{-6}$ & $2.5\cdot10^{-7}$
& $1\cdot10^{-8}$ & $1\cdot10^{-8}$ \\ \hline
$N$ & 34009 & 32980 & 36655 & 33421 & 34450 & 30040 \\ \hline
$\lambda_1$ &$354_{43.169162}^{24.194986}$ & $158_{99.470006}^{84.340834}$ & 
$89_{14.1217006}^{00.8711417}$ & $57_{28.5030533}^{16.2444537}$ &
$40_{33.0287132}^{21.2855619}$ & $30_{29.4734311}^{17.9563370}$ \\ \hline \hline
$s_1$ & $0.3097_{337281}^{751958}$ & $0.2310_{241412}^{791320}$ & $0.170_{6976522}^
{7611453}$ & $0.123_{6657365}^{7319839}$ & $0.086_{8374821}^{9008096}$
& $0.0582_{044100}^{598605}$ \\ \hline
$r_1$ & $0.2792_{012640}^{386428}$ & $0.2021_{168919}^{650011}$ & 
$0.145_{1659786}^{2199751}$ & $0.1024_{057897}^{606481}$ & 
$0.0701_{195823}^{707180}$ & $0.0459_{009493}^{446784}$ \\ \hline \hline
$s_2$ & $0.1459_{659148}^{854569}$ & $0.0835_{750755}^{949688}$ & 
$0.046_{4973316}^{5146269}$ & $0.0247_{301082}^{433560}$ & 
$0.0123_{155385}^{245198}$ & $0.0055_{776108}^{829245}$ \\ \hline
$r_2$ & $0.1316_{369362}^{545587}$ & $0.0731_{289927}^{463995}$ & 
$0.0395_{463421}^{610519}$ & $0.0204_{778925}^{888625}$ & 
$0.0099_{450264}^{522789}$ & $0.004_{3982853}^{4024756}$ \\ \hline \hline
$s_3$ & $0.0691_{063846}^{156363}$ & $0.0302_{888290}^{960386}$ & 
$0.01267_{35136}^{82277}$ & $0.00494_{62962}^{89459}$ & 
$0.00174_{67069}^{79807}$ & $0.00053_{44933}^{50025}$ \\ \hline
$r_3$ & $0.0623_{221671}^{300142}$ & $0.02650_{20383}^{83465}$ &
$0.0107_{792514}^{832608}$ & $0.00409_{59277}^{81219}$ & 
$0.00141_{04601}^{14887}$ & $0.000421_{5086}^{9102}$ \\ \hline\hline
$s_4$ & $0.03272_{54060}^{97872}$ & $0.0109_{77458}^{800588}$ & 
$0.00345_{43607}^{56456}$ & $0.000989_{3135}^{8434}$ & 
$0.000247_{7346}^{9152}$ & $0.0000512_{196}^{684}$ \\ \hline
$r_4$ & $0.02951_{38412}^{77925}$ & $0.00960_{54332}^{77195}$ & 
$0.00293_{80625}^{91553}$ & $0.000819_{2651}^{7040}$ &
$0.000200_{0575}^{2034}$ & $0.000040_{3912}^{4297}$ \\ \hline\hline
\end{tabular}}

\medskip

{\small\begin{tabular}{|c|c|c|c|c|c|c|} \hline
$\theta$  & $90^\circ$ & $100^\circ$ & $110^\circ$ & $120^\circ$ &
$130^\circ$ & $140^\circ$ \\ \hline
$h_1$ & 0.199 & 0.222 & 0.245 & 0.268 & 0.291 & 0.315 \\ \hline
$\rho_1$ &  $1\cdot10^{-8}$ & $5\cdot10^{-9}$ &
$2.5\cdot10^{-9}$ & $5\cdot10^{-10}$ & $5\cdot10^{-10}$ & $5\cdot10^{-10}$ \\ 
\hline
$N$ &  29746 & 24601 & 23092 & 22651 & 20887 & 19270 \\ \hline
$\lambda_1$ &  $23_{87.7996658}^{76.3155992}$
& $19_{52.8971428}^{41.3070303}$ & $16_{44.4979497}^{32.6947173}$ &
$14_{17.8265343}^{05.7227746}$ & $12_{46.3582958}^{33.8788513}$ &
$11_{13.5880500}^{00.6652772}$ \\ \hline \hline
$s_1$ &  $0.0364_{024346}^{463359}$
& $0.0204_{616208}^{920930}$ & $0.0096_{399984}^{573740}$ & $0.0032_{667832}^
{737930}$ & $0.00051_{28490}^{41408}$ & $\!0.0000035_{312}^{415}\!$ \\ \hline
$r_1$ &  $0.028_{0724689}^{1063243}$ & $0.0154_{513578}
^{743684}$ & $0.0071_{369993}^{498633}$ & $0.00237_{42004}^{92947}$ &
$0.00036_{62871}^{72098}$ & $\!0.00000248_{15}^{88}\!$ \\ \hline \hline
$s_2$ &  $0.00219_{72279}^{98778}$ & $0.000_{6989686}
^{7000095}$ & $0.000156_{2557}^{5374}$ & $0.000018_{0886}^{1274}$
& $\!0.00000045_{01}^{12}\!$ & --- \\ \hline
$r_2$ &  $0.00169_{44661}^{65096}$ & $0.00052_{78278}
^{86069}$ & $0.000115_{6862}^{8947}$ & $0.0000131_{461}^{743}$
& $\!0.00000032_{15}^{23}\!$ & --- \\ \hline \hline
$s_3$ & $0.000132_{6235}^{7834}$ & $0.000023_{8767}
^{9123}$ & $0.00000253_{28}^{73}$ & $\!0.000000100_{2}^{4}\!$
& --- & --- \\ \hline
$r_3$ & $0.000102_{2789}^{4022}$ & $0.0000180_{306}
^{567}$ & $0.00000187_{52}^{86}$ & $\!0.000000072_{8}^{9}\!$ &
--- & --- \\ \hline\hline
$s_4$ & $0.0000080_{051}^{147}$ & $\!0.00000081_{56}
^{68}\!$ & $\!0.000000041_{0}^{1}\!$ & --- & --- & --- \\ \hline
$r_4$ & $0.0000061_{736}^{810}$ & $\!0.00000061_{59}
^{68}\!$ & $\!0.000000030_{4}^{5}\!$ & --- & --- & --- \\ \hline\hline
\end{tabular}}
\caption{Local zeros ($s_n$) and local extrema ($r_n$), as calculated
on the coarse meshes, of the first eigenvalue for the clamped plate
eigenvalue problem on a circular sector of angle $\theta$. \label{t:3.1}}
\end{center}
\end{table}
\begin{table}[htb]
\begin{center}
{\small\begin{tabular}{|c|c|c|c|c|c|c|c|} \hline
$\alpha$ & $30^\circ$ & $40^\circ$ & $50^\circ$ &
$60^\circ$ & $70^\circ$ & $80^\circ$ \\ \hline
$h_1$ & 0.033 & 0.044 & 0.055 & 0.066 & 0.077 & 0.087 \\ \hline
$\rho_1$ &  $1\cdot10^{-2}$ & $2\cdot10^{-3}$ & $2.5\cdot10^{-4}$ & 
$5\cdot10^{-5}$ & $1\cdot10^{-5}$ & $1\cdot10^{-6}$ \\ \hline
$N$ & 82084 & 82627 & 87865 & 87283 & 86701 & 90775 \\ \hline
$\lambda_1$ & $354_{33.606980}^{28.863922}$ & $158_{91.837688}^{88.056086}$ & 
$890_{7.4275436}^{4.1158527}$ & $57_{22.2999555}^{19.2365745}$ &
$402_{7.0760866}^{4.1419602}$ & $302_{3.6247696}^{0.7476335}$ \\ \hline \hline
$s_1$ & $0.3097_{547154}^{650827}$ & $0.2310_{519414}^{656892}$ & $0.1707_{293299}^
{452024}$ & $0.123_{6992733}^{7158333}$ & $0.0868_{694257}^{852561}$
& $0.0582_{326311}^{464922}$ \\ \hline
$r_1$ & $0.2792_{262612}^{356059}$ & $0.2021_{346638}^{466898}$ &
$0.145_{1981906}^{2116893}$ & $0.1024_{314713}^{451849}$ & 
$0.0701_{473169}^{601000}$ & $0.0459_{221772}^{331081}$ \\ \hline \hline
$s_2$ & $0.1459_{757551}^{806409}$ & $0.08358_{49350}^{99082}$ & 
$0.0465_{059576}^{102813}$ & $0.0247_{367298}^{400416}$ & 
$0.01232_{00774}^{23226}$ & $0.0055_{792134}^{801286}$ \\ \hline
$r_2$ & $0.13164_{45423}^{89477}$ & $0.0731_{370529}^{414045}$ 
& $0.03955_{38988}^{75761}$ & $0.02048_{42106}^{69531}$ & 
$0.0099_{485408}^{503538}$ & $0.00440_{05850}^{16326}$ \\ \hline \hline
$s_3$ & $0.06911_{09631}^{32761}$ & $0.03029_{24195}^{42219}$ & 
$0.01267_{58803}^{70587}$ & $0.00494_{76260}^{82884}$ & 
$0.001747_{3480}^{6664}$ & $0.000534_{7501}^{8774}$ \\ \hline
$r_3$ & $0.06232_{69688}^{90547}$ & $0.02650_{60380}^{76150}$ & 
$0.01078_{09985}^{20009}$ & $0.004097_{0420}^{5905}$ & 
$0.00141_{09998}^{12570}$ & $0.000421_{7013}^{8016}$ \\ \hline\hline
$s_4$ & $0.03272_{77538}^{88491}$ & $0.01097_{87694}^{94226}$ & 
$0.003455_{0087}^{3299}$ & $0.000989_{5811}^{7136}$ & 
$0.0002478_{252}^{703}$ & $0.0000512_{442}^{564}$ \\ \hline
$r_4$ & $0.02951_{52608}^{62568}$ & $0.00960_{64861}^{70576}$ &
$0.002938_{5291}^{8023}$ & $0.000819_{4572}^{5669}$ & 
$0.0002001_{219}^{584}$ & $0.0000404_{110}^{206}$ \\ \hline
\end{tabular}}

\medskip

{\small\begin{tabular}{|c|c|c|c|c|c|c|} \hline
$\theta$ & $90^\circ$ & $100^\circ$ & $110^\circ$ & $120^\circ$ &
$130^\circ$ & $140^\circ$ \\ \hline
$h_1$ & 0.098 & 0.110 & 0.121 & 0.132 & 0.143 & 0.154 \\ \hline
$\rho_1$ &  $2.5\cdot10^{-7}$ & $5\cdot10^{-8}$ &
$1\cdot10^{-8}$ & $2.5\cdot10^{-9}$ & $5\cdot10^{-10}$ & $5\cdot10^{-10}$ \\ 
\hline
$N$  & 88447 & 87865 & 87283 & 85828 & 85246 & 79096 \\ \hline
$\lambda_1$ &  $23_{81.9568311}^{79.0885238}$
& $194_{6.9890139}^{4.0948765}$ & $163_{8.4693656}^{5.5227545}$ &
$14_{11.6320780}^{08.6112832}$ & $123_{9.9585482}^{6.8449411}$ &
$110_{6.9471887}^{3.7240407}$ \\ \hline \hline
$s_1$ &  $0.0364_{247712}^{357449}$
& $0.0204_{770656}^{846824}$ & $0.0096_{488600}^{532030}$ & $0.00327_{03507}^
{21026}$ & $0.000513_{5089}^{8318}$ & $\!0.00000353_{65}^{91}\!$ \\ \hline
$r_1$ &  $0.02809_{01923}^{86552}$ & $0.01546_{30180}
^{87697}$ & $0.00714_{36744}^{68898}$ & $0.00237_{67308}^{80040}$ &
$0.000366_{7660}^{9966}$ & $\!0.00000248_{53}^{71}\!$ \\ \hline \hline
$s_2$ & $0.00219_{85741}^{92365}$ & $0.000699_{4954}
^{7556}$ & $0.000156_{3990}^{4694}$ & $0.0000181_{083}^{180}$
& $\!0.000000450_{7}^{9}\!$ & --- \\ \hline
$r_2$ &  $0.00169_{55059}^{60167}$ & $0.000528_{2165}
^{4129}$ & $0.000115_{7914}^{8435}$ & $0.00001316_{03}^{74}$
& $\!0.00000032_{19}^{21}\!$ & --- \\ \hline \hline
$s_3$ & $0.0001327_{046}^{446}$ & $0.000023_{8947}
^{9036}$ & $0.00000253_{51}^{62}$ & $\!0.000000100_{2}^{3}\!$
& --- & --- \\ \hline
$r_3$ & $0.0001023_{391}^{699}$ & $0.0000180_{439}
^{506}$ & $0.00000187_{69}^{73}$ & $\!0.00000007_{29}^{30}\!$ &
--- & --- \\ \hline\hline
$s_4$ &  $0.00000801_{00}^{24}$ & $\!0.000000816_{2}
^{5}\!$ & $\!0.0000000411\!$ & --- & --- & --- \\ \hline
$r_4$ &  $0.00000617_{71}^{90}$ & $\!0.000000616_{4}
^{5}\!$ & $\!0.0000000304\!$ & --- & --- & --- \\ \hline 
\end{tabular}}
\caption{Local zeros ($s_n$) and local extrema ($r_n$), as calculated
on the fine meshes,  of the first
eigenvalue  for the clamped plate eigenvalue problem on a circular
sector of angle $\theta$. \label{t:3.2}}
\end{center}
\end{table}

\noindent In each case
the numerical results quoted take the form of pairs of values corresponding
to those calculated on the inner (bottom value) and outer (top value)
polygons respectively.
It should be noted that, with one exception, each pair of values
in Table \ref{t:3.2} lies between the corresponding values in Table \ref{t:3.1}.
While this proves nothing rigorous about bounds on the eigenvalue or the positions
of the zeros and the extrema, it does provide a degree of
confidence in the calculations. Furthermore, the one value
($r_3$ when $\theta\!=\!120^\circ$) for which there is any discrepancy between the 
two sets of results is clearly on the limit of our computational accuracy. Even
in this case however, all four estimates of this value agree to two significant
figures and have an absolute difference of less than $ 10^{-9}$.
The results in Tables \ref{t:3.1} and \ref{t:3.2} clearly show that,
as the theory outlined above predicts, the frequency of the oscillations
of the eigenfunction, $u$, decreases as $\theta$ increases.
Table \ref{t:3.3} further illustrates this point by presenting
the ratios of the positions of consecutive zeros and extrema for
the fine mesh calculations. This table also shows the known asymptotic
limit for these ratios in each case. It is notable that these asymptotic ratios
are remarkably well approximated even by the first few ratios that we are able
to calculate here.
Furthermore, as $\theta$ approaches an angle of $140^\circ$, it is clear that
our numerical calculations are having great difficulty in resolving even the
second change of sign. This is to be expected since the solution of equation
(\ref{eq:1.4}) with smallest real part is equal to $2.826869+0.261695 i$ for this 
value
of $\theta$: hence, by (\ref{eq:1.5}) and (\ref{eq:1.6}), the limiting ratio of
the position and magnitude of consecutive
extrema is equal to  $163533.23$ and $5.472 \times 10^{14}$ respectively.
Given that $r_1$ is already very small and
$|t_1|\!\approx\!1.30729\cdot 10^{-15}$ times the maximum magnitude of 
$u$,
double precision calculations could not reasonably be expected
to locate the second extremum with any significant figures of accuracy.
For this reason, in Table \ref{t:3.4} we are only able to present calculations
of the location of the first zero of $u$ when considering values of
$\theta$ greater than $140^\circ$. The parameters used for the generation
of the meshes for these calculations are $\rho_1\!=\!5 \cdot10^{-10}$ and 
$h_1\!\approx\!0.16$,
which, in principal at least, allow us to observe a sign change at a distance of
$O(10^{-11})$ from the corner. In practice however we are unable to detect
the first change in sign when $\theta\!=\!145^\circ$ and $\theta\!=\!146^\circ$
using such meshes. Nevertheless it is clear that these simple calculations
do verify the theory which predicts the loss of oscillations at
some critical angle $\theta_c$, \mbox{even though we are unable to accurately
determine a value for $\theta_c$ using these meshes.}

\begin{table}[htb]
\begin{center}
{\small \begin{tabular}{|c|c|c|c|c|c|c|}
\hline $\theta$ & $30^\circ$ & $40^\circ$ & $50^\circ$ &
$60^\circ$ & $70^\circ$ & $80^\circ$  \\ \hline
$s_1/s_2$ &  2.1219600 & 2.7642773 & 3.6711281 & 5.0006316 &
7.0510454 & 10.435415 \\ \hline
$s_2/s_3$ &  2.1121939 & 2.7592690 & 3.6688543 & 4.9997170 & 7.0507290 
& 10.435321 \\ \hline
$s_3/s_4$ & 2.1116928 & 2.7591817 & 3.6688418 & 4.9997176 & 7.0507283 & 
10.435323 \\ \hline \hline
limit & 2.1116888 & 2.7591817 & 3.6688419 & 4.9997171 & 7.0507288 &
10.435321 \\ \hline\hline
$r_1/r_2$ & 2.1210622 & 2.7637792 & 3.6708945 & 5.0005086 &
7.0510156 & 10.435471 \\ \hline
$r_2/r_3$ &  2.1121602 & 2.7592601 & 3.6688530 & 4.9997561 & 7.0507030 
& 10.435314 \\ \hline
$r_3/r_4$ & $2.11168_{62}^{56}$ & 2.7591814 & 3.6688418 & 4.9997024 
& 7.0507288 & 10.435311 \\ \hline \hline
limit & 2.1116888 & 2.7591817 & 3.6688419 & 4.9997171 & 7.0507288 &
10.435321 \\ \hline
\end{tabular}

\medskip

\begin{tabular}{|c|c|c|c|c|c|c|}
\hline $\theta$ & $90^\circ$ & $100^\circ$ & $110^\circ$ &
$120^\circ$ & $130^\circ$ & $140^\circ$  \\ \hline
$s_1/s_2$ &  16.567452 & 29.274052 & 61.693788 & 180.59907 &
\,$1139.4_{724}^{505}$\, & \ \ \ \,\,\,---\ \ \ \,\,\,\\ \hline
$s_2/s_3$ &  16.567430 & 29.274040 & 61.693978 & $180.5_{9586}^{8995}$
& --- & --- \\ \hline
$s_3/s_4$ &  $16.567425$ & $29.274_{153}^{355}$ & $61._{673203}^
{718979}$ & --- & --- & --- \\ \hline \hline
limit & 16.56743 &29.27404& 61.69387& 180.5992& 1139.464& 163533.2 \\ 
\hline
\hline
$r_1/r_2$ &  16.567440 & 29.274017 & 61.694352 & $180.598_{34}^{28}$ &
\,$1139.4_{450}^{312}$\, & \ \ \ \,\,\,---\ \ \ \,\,\,\\ \hline
$r_2/r_3$ &  16.567535 & $29.27394_0^8$ & $61.693_{744}^{540}$ &
$180.62_{716}^{481}$ & --- & --- \\ \hline
$r_3/r_4$ &  $16.5674_{25}^{17}$ & $29.274_{105}^{727}$ &
$61.7_{51464}^{30719}$ & --- & --- & --- \\ \hline \hline
limit & 16.56743 &29.27404& 61.69387& 180.5992& 1139.464& 163533.2 \\ 
\hline
\end{tabular}}
\caption{Ratios of consecutive local zeros and extrema of the eigenfunction
(along with the theoretical asymptotic value)
for different values of $\theta$. \label{t:3.3} }
\end{center}
\end{table}
\begin{table}[htb]
\begin{center}
{\footnotesize\begin{tabular}{|c|c|c|c|c|c|c|c|} \hline
$\theta$ & $141^\circ$ & $142^\circ$ & $143^\circ$ & $144^\circ$ & $145^\circ$ &
$146^\circ$ & $147^\circ$ \\ \hline
$s_1$ & $0.00000111_{165}^{248}$ & $0.000000240_{47}^{66}$
& $0.0000000275_{0}^{1}$ & $0.000000000_{85}^{93}$ & \mbox{\qquad ---\qquad} 
& \mbox{\qquad ---\qquad} & \mbox{\qquad ---\qquad} \\ \hline
\end{tabular}}
\caption{Calculated positions of  the  local zeros of the eigenfunction
for values of $\theta$ above $140^\circ$. \label{t:3.4} }
\end{center}
\end{table}

Another interesting situation to investigate numerically would be the
limiting behaviour of the principal eigenfunction as $\theta \rightarrow 0$.
Theory predicts that, as $\theta$ decreases, oscillations will become
more frequent and their amplitude will decay more quickly as a function of
$r$ (although the asymptotic ratio of successive eigenvalues will decrease).
This presents different, but equally demanding, computational challenges
to those for $\theta \rightarrow \theta_c$. In particular it is no longer
sufficient to refine the mesh heavily only in the neighbourhood
of the corner since changes in sign occur frequently throughout a large
part of the domain. Hence, for very small angles $\theta$, a very dense
finite element mesh is required over almost all of $\Omega_\theta$ which,
despite the narrowing of the domain,
leads to very large algebraic systems (\ref{e2.7}).
In fact, even for $\theta=10^\circ$, we observe a reduction in the
apparent rate of convergence when considering the position of
the zeros of $u$ on sequences of finer meshes.
This leads us to believe that a thorough numerical investigation
of the behaviour of this eigenfunction as $\theta \rightarrow 0$
will be an even more computationally demanding task.

In the numerical results that we present above it has been possible,
for the first time that we are aware, to verify a number of
theoretical asymptotic results concerning the existence of nodal lines
in the vicinity of a corner  with arbitrary internal angle $\theta$.
Furthermore we have computed, with some confidence,
estimates for the exact locations of zeros and extrema of the principal 
eigenfunction
for a number of different values of $\theta$.
More comprehensive results, including those for the buckling plate
problem, (\ref{eq:1.2}), can also be found in \cite{MilanT}.
In the following section we move on from providing numerical verification
and extensions of known analytic results to demonstrate that numerics
can also give some insight into problems for which there are
no theoretical predictions.

\section{Parity property of eigenfunctions} \label{s:4}
In this section we again consider problems (\ref{eq:1.1}) and (\ref{eq:1.2})
but now on another family of domains, $\Omega_c$ say, which are defined
as the regions bounded by the curves (\ref{eq:1.7}) for $c >0$.
We are particularly interested in the behaviour of the eigenfunctions
which correspond to the smallest eigenvalues as $c \rightarrow 0$.
In this limit the connected domain $\Omega_c$ tends to a disconnected
domain which is symmetric about the $y-$axis. It follows therefore
that each eigenvalue must have an even multiplicity and, in particular,
the smallest eigenvalue will be repeated. When $c>0$, but small, therefore,
we expect the eigenvalues of the biharmonic operator to
occur in distinct pairs (i.e.\ there will be a small gap between the
$(2k-1)^{th}$ and the $2k^{th}$ eigenvalues for positive integers $k$).
The same argument holds when considering the {\em Laplacian} eigenvalue problem
on $\Omega_c$, subject to zero Dirichlet boundary conditions.
Indeed, it is known that for this second order problem the
eigenfunction corresponding to the lower eigenvalue in each pair
is always an even function of $x$ while the other eigenfunction
in the pair is always odd. No such result is known in the case of the
biharmonic operator however. In Table \ref{t:4.1} below  we present 
numerical evidence that such a result is false. This is achieved by 
demonstrating that branches of the two least eigenvalues (corresponding
to eigenfunctions of different parity), viewed as functions of $c$, 
appear to cross on more than one occasion as $c$ approaches zero.

\begin{table}[htb]
\begin{center}
{\small\begin{tabular}{ | c | c | c | c | c |} \hline
$c$ & $ \lambda_{even} / \lambda_{odd}$ & $\lambda_{even}$ & $\lambda_{odd}$ & 
$N$\\ \hline \hline
1.0 & 0.304005 & 51.9393395 & 170.850197 & 129602\\ \hline
0.9 & 0.349358 & 68.2814038 & 195.448422 & 141122\\ \hline
0.8 & 0.410835 & 94.5437752 & 230.125842 & 152642 \\ \hline
0.7 & 0.494593 & 139.440088 & 281.928749 & 108290 \\ \hline
0.6 & 0.608020 & 221.940406 & 365.021620 & 119810 \\ \hline
0.5 & 0.755019 & 385.377821 & 510.421666 & 133634 \\ \hline
0.4 & 0.914099 & 724.064578 & 792.107738 & 112754 \\ \hline
0.35& 0.972477 & 1005.87754 & 1034.34630 & 120530 \\ \hline
0.325&0.990509 & 1185.48005 & 1196.83874 & 128306 \\ \hline
0.3 & 1.000849& 1397.46767 & 1396.29291 & 134138 \\ \hline
0.275&1.004906& 1650.06652 & 1642.01017 & 141914 \\ \hline
0.25 &1.004943& 1955.47056 & 1945.85190 & 151634 \\ \hline
0.225&1.003236& 2330.95499 & 2323.43715 & 161354 \\ \hline
0.2 & 1.001465& 2800.48607 & 2796.38901 & 136514 \\ \hline
0.175&1.000403& 3396.70646 & 3395.33941 & 125138 \\ \hline
0.15 &1.000031& 4165.22544 & 4165.09438 & 136226 \\ \hline
0.125&0.999987& 5172.36322 & 5172.42820 & 152066 \\ \hline
0.1  &0.999998& 6518.06051 & 6518.07172 & 144002 \\ \hline
0.075&1.000000& 8356.71703 & 8356.71678 & 136082 \\ \hline
0.05 &1.000000& 10936.0178 & 10936.0178 & 133634 \\ \hline
0.025&1.000000& 14670.7338 & 14670.7338 & 126362 \\ \hline
0.01 &1.000000& 17752.3432 & 17752.3432 & 118082 \\ \hline
\end{tabular}}
\caption{Estimates of the two smallest eigenvalues, and their parity, for
the clamped plate problem on $\Omega_c$. \label{t:4.1}}
\end{center}
\end{table}

The results presented have been calculated using the same
piecewise cubic  mixed finite element method described and used in the previous
sections. In practice our computations were performed on a discretization
of $\Omega_c$ in which the curved boundary is represented
by an interpolating polygon. Since our motivation here is primarily to find
an example in which the parity of the eigenfunctions changes with the
domain geometry, this approximation is of no significant consequence.
We again make use of unstructured meshes of triangles in the interior
of the domain and use $N$ to represent the total number of degrees
of freedom in the finite element approximation. For each value of $c$
calculations were performed on a sequence of meshes with different
numbers of triangles. Only the results of computations on the finest
meshes are presented in Table \ref{t:4.1}, however an indication of
the reliability of these figures is given by comparison with
the corresponding calculations on coarser meshes. In particular,
Table~\ref{t:4.2} gives the relative differences (as percentages)
between consecutive ratios $\lambda_{even} / \lambda_{odd}$ computed on
four different meshes for each value of the parameter $c$. In the table
$N_i$ ($i=1,\ldots,4$) denotes the total dimension of the linear system
arising from the four different discretizations of the domain, and
$\delta\,(\%)$ stands for the relative difference in the ratios
$\lambda_{even} / \lambda_{odd}$ computed on consecutive meshes.
It is reasonable to expect that our results are reliable if this
relative ``error'' decreases (and becomes sufficiently small) as the mesh
is refined. As we can see from Table~\ref{t:4.2} this is
the case for all meshes and all values of the parameter $c$ considered.

\begin{table}[htb]
\begin{center}
{\footnotesize\begin{tabular}{|c|c|c|c|c|c|c|c|c|c|c|c|} \hline
$c$ & 1.0 & 0.9 & 0.8 & 0.7 & 0.6 & 0.5 & 0.4 & 0.35 & 0.325 & 0.3 & 0.275 \\ 
\hline\hline
$N_1$ & 52202 & 54002 & 59402 & 41762 & 47522 & 53282 & 50546 & 54434 & 57026
& 60914 & 63506 \\ \hline
$N_2$ & 73442 & 79922 & 86402 & 60482 & 67394 & 76034 & 68042 & 74090 & 78626
& 81650 & 86186 \\ \hline
$\delta\,(\%)$ & 1.18e-2 & 1.46e-2 & 1.29e-2 & 1.97e-2 & 1.45e-2 & 1.21e-2 
& 6.15e-3 & 3.24e-3 & 1.75e-3 & 7.70e-4 & 2.47e-4 \\ \hline\hline 
$N_2$ & 73442 & 79922 & 86402 & 60482 & 67394 & 76034 & 68042 & 74090 & 78626
& 81650 & 86186 \\ \hline
$N_3$ & 98282 & 105842 & 115922 & 80642 & 90722 & 102818 & 89858 & 95042 
& 101954 & 105410 & 112322 \\ \hline
$\delta\,(\%)$ & 7.47e-3 & 6.81e-3 & 7.77e-3 & 9.19e-3 & 8.95e-3 & 6.14e-3 
& 4.84e-3 & 1.79e-3 & 1.16e-3 & 5.19e-4 & 1.38e-4 \\ \hline\hline 
$N_3$ & 98282 & 105842 & 115922 & 80642 & 90722 & 102818 & 89858 & 95042 
& 101954 & 105410 & 112322 \\ \hline
$N_4$ & 129602 & 141122 & 152642 & 108290 & 119810 & 133634 & 112754 & 120530
& 128306 & 134138 & 141914 \\ \hline
$\delta\,(\%)$ & 5.03e-3 & 5.89e-3 & 5.05e-3 & 7.93e-3 & 5.90e-3 & 5.01e-3 
& 2.69e-3 & 1.53e-3 & 8.05e-4 & 4.19e-4 & 1.16e-4 \\ \hline 
\end{tabular}}

\medskip

{\footnotesize\begin{tabular}{|c|c|c|c|c|c|c|c|c|c|c|c|} \hline
$c$ & 0.25 & 0.225 & 0.2 & 0.175 & 0.15 & 0.125 & 0.1 & 0.075 & 0.05 & 0.025 & 0.01\\ 
\hline\hline
$N_1$ & 68690 & 72758 & 76466 & 66818 & 73730 & 80642 & 70562 & 82658 & 52562
& 74882 & 58466 \\ \hline
$N_2$ & 92234 & 96770 & 95042 & 84242 & 93314 & 102386 & 92162 & 108290 & 75170
& 91082 & 76034 \\ \hline
$\delta\,(\%)$ & 3.02e-5 & 8.81e-5 & 4.88e-5 & 2.42e-5 & 4.80e-6 & 9.99e-8 
& 9.99e-8 & 0.000 & 0.000 & 0.000 & 0.000 \\ \hline\hline 
$N_2$ & 92234 & 96770 & 95042 & 84242 & 93314 & 102386 & 92162 & 108290 & 75170
& 91082 & 76034 \\ \hline
$N_3$ & 119234 & 127874 & 115634 & 103682 & 112322 & 126722 & 116642 & 122402 
& 101810 & 108002 & 95906 \\ \hline
$\delta\,(\%)$ & 2.01e-5 & 6.43e-5 & 3.64e-5 & 1.79e-5 & 3.10e-6 & 9.99e-8 
& 9.99e-8 & 0.000 & 0.000 & 0.000 & 0.000 \\ \hline\hline 
$N_3$ & 119234 & 127874 & 115634 & 103682 & 112322 & 126722 & 116642 & 122402 
& 101810 & 108002 & 95906 \\ \hline
$N_4$ & 151634 & 161354 & 136514 & 125138 & 136226 & 152066 & 144002 & 136082 
& 133634 & 126362 & 118082 \\ \hline
$\delta\,(\%)$ & 1.41e-5 & 4.42e-5 & 2.78e-5 & 1.50e-5 & 2.60e-6 & 0.000 
& 0.000 & 0.000 & 0.000 & 0.000 & 0.000 \\ \hline 
\end{tabular}}
\caption{Percentage discrepancies $\delta$ between the ratios
$\lambda_{even}/\lambda_{odd}$ obtained on four different meshes for each 
domain $\Omega_c$ (discrepancies of less than 1.0e-8 have been rounded
to zero). \label{t:4.2}}
\end{center}
\end{table}

Table \ref{t:4.1} clearly shows that, for our finite element
discretization at least,  there are indeed values of $c$ for which
the branches of the two smallest eigenvalues cross. Inspection of
the corresponding eigenfunctions demonstrates that one of these functions
is always even and the other is always odd: we refer to the
associated eigenvalues as $\lambda_{even}(c)$ and $\lambda_{odd}(c)$
respectively. Further numerical computations on a sequence of coarser
grids, \cite{MilanT}, suggest that as the mesh size $h$ is reduced,
the locations of these crossing points remain stable and the
difference between $\lambda_{even}(c)$ and $\lambda_{odd}(c)$
does not tend to zero for values of $c$ between these points.
This suggests that it is reasonably safe for us to conclude that the
parity of the eigenfunction corresponding to the least eigenvalue
of (\ref{eq:1.1}) does indeed change as the domain $\Omega_c$ evolves.
Although we have only presented results for the clamped plate problem
here, it is shown in \cite{MilanT} that similar behaviour may be
observed for the buckling plate problem (\ref{eq:1.2}) on these domains.
It is interesting to hypothesize on the behaviour of the eigenvalue branches
as $c$ gets very close to its limiting value of zero. The resolution of the
computations that we have undertaken only appears to allow us to identify
the first two or three values of $c$ for which the least eigenvalue is
repeated.  Beyond this point discretization and rounding errors make it
impossible to distinguish between them, even though they are almost
certainly different  for most values of $c>0$. There may however be
further values of $c$, smaller than those that we have been able to
identify, for which $\lambda_{even}(c)=\lambda_{odd}(c)$: we conjecture
that there are an infinite number of such values.

\section{Discussion} \label{s:5}
In this paper our aim has been to demonstrate that the finite element method
may be successfully applied to biharmonic eigenvalue problems of the form 
(\ref{eq:1.1}) and (\ref{eq:1.2}), and that it may be used both to verify
existing asymptotic results and to allow new conjectures to be made about
the behaviour of their eigenfunctions. An outline of the finite element
implementation is given in Section \ref{s:2}, with complete details
available in \cite{BJM99}. Verification of the accuracy of this approach
is provided by comparison with known analytical results applied on a unit
square domain, along with the recent high precision numerical results
(calculated using an entirely different discretization scheme) that appear
in \cite{Bjo}.

An advantage of the mixed finite element method used in this work is that
it may easily be applied on unstructured grids of triangles, ${\cal{T}}^h$,
and therefore on arbitrary polygonal domains in $\Re^2$. We have made use
of this fact to investigate the spectral properties of the biharmonic
operator on two further geometries: sectors of the unit circle and a family
of non-convex domains. Furthermore it is possible to generate meshes for
which the degrees of freedom in the corresponding piecewise cubic trial
spaces are not uniformly distributed throughout $\Omega$. Where appropriate,
this allows us to obtain a greater resolution of the eigenfunctions in
those regions in which we wish to focus, without affecting either the
dimension or the sparsity of the discrete linear system (\ref{e2.7}).

The results presented in Section \ref{s:3} are in agreement with the
asymptotic theory that appears in
\cite{B+R},\cite{Coff},\cite{Kond},\cite{KKM} and also provide quantitative
estimates of the positions of a number of sign changes and extremal
points for various values of the internal angle $\theta$. In Section
\ref{s:4} we extend the numerical investigation further by considering
a family of problems for which the eigenfunctions have significantly
different qualitative properties. This has led to evidence that, unlike
the Laplace operator, the parity of the least eigenfunction of the
biharmonic operator can change as the domain is deformed.

The unstructured finite element approach that we have used in this
work may also be applied to further problems in this area. Clearly it
is possible to investigate the spectral properties of the biharmonic
operator on a greater variety of geometries in 2-d, and subject to a
wider variety of boundary conditions. This would include, for example,
further studies in the circular sector domains for very small angles
$\theta$. Perhaps more interestingly however (and certainly more
computationally demanding) would be extensions to higher order
operators and/or larger numbers of dimensions.

We conclude this paper by noting that none of the computational results
that are presented here provide {\em rigorous} information about the
eigenfunctions considered (although the apparent mesh convergence
observed is certainly a strong indicator). It would be interesting
therefore to attempt to obtain provably correct bounds for certain
quantities associated with the eigenvalues and eigenfunctions of the
biharmonic operator on the regions considered. In \cite{Wien}, for
example, Weiners has made a start on this problem by computing
enclosures for the first eigenfunction in the neighbourhood of a corner
of the unit square. This involves computing estimates for the Sobolev
embedding constants which in turn requires an upper bound for the defect of
the eigenfunction. He is able to show that the defect is less than
$0.000000893$ which is sufficiently accurate to allow him to conclude that
the first eigenfunction changes sign. Since the estimates of the amplitude
of the second and higher sign changes that we have found are
considerably smaller than for the first, new methods are likely to be
required in order to compute the defect with sufficient precision.

 \section*{Acknowledgements}
BMB and MDM acknowledge the
EPRSC for support under grant GR/K84745 and PKJ acknowledges the support of the
Leverhulme Trust.

\end{document}